\documentclass[oneside]{amsart}

  \newcommand{\Section}{\section}
  \newcommand{\SubSection}{\subsection}

  \newcommand{\rar}{\rightarrow}

  \newcommand{\scr}[1]{\mathscr{#1}}

\newcommand{\bb}[1]{\mathbb{#1}}
\newcommand{\im}[0]{ \mathbf{i} }

\usepackage{amssymb,amsfonts,mathrsfs}
\usepackage[all,arc]{xy}
\usepackage{enumerate}
\usepackage{mathrsfs}

\usepackage{comment}

\usepackage[letterpaper]{geometry}  

\usepackage{enumitem}   

\usepackage[utf8]{inputenc}

\usepackage{graphicx} 

\usepackage{hyperref}
\hypersetup{
    colorlinks,
    citecolor=black,
    filecolor=black,
    linkcolor=black,
    urlcolor=black
}

\theoremstyle{definition} 
\newtheorem{thm}{Theorem}[section]

\newtheorem{lem}[thm]{Lemma}

\newtheorem{fact}[thm]{Fact}

\theoremstyle{definition}
\newtheorem{defin}[thm]{Definition}

\newtheorem{notn}[thm]{Notation} 
\newtheorem{conv}[thm]{Convention} 

\theoremstyle{remark}
\newtheorem{rem}[thm]{Remark}

\makeatletter
\let\c@equation\c@thm
\makeatother
\numberwithin{equation}{section}

\title{ Rigidity via modular properties of theta functions } 

\author{Indraneel Tambe}

\date{}

\begin{document}

\maketitle


In this paper we use methods of Liu
    \cite{KL-Mod, KL-EGTF}
to show that the twisted Dirac operators $\mathscr{D}$
on certain bundles $\Phi$
considered by Guan and Wang
\cite{GuanWangII} 
are rigid. 
    To do so, 
    we use a Lefschetz formula and Atiyah-Bott localization 
to obtain formulas for the Lefschetz numbers $\mathscr{L}$ of these operators $\mathscr{D}$
in terms of Jacobi theta functions;
  then, using the translational and modular transformation properties of theta functions and the properties of their zeros,
we prove $\mathscr{L}$ is constant 
    provided certain conditions on characteristic classes hold,
thus showing the rigidity of $\mathscr{D}$ on $\Phi$ under these conditions.


\tableofcontents


\Section{Introduction}

Let $G$ be a Lie group and let $P$ be a $G$-equivariant Fredholm differential operator acting on the sections of a $G$-equivariant smooth vector bundle $W \rightarrow M$.
The kernel and cokernel of $P$ are finite-dimensional representations of $G$,
and we obtain a well-defined function $L$ on $G$ given by $L(g) := \mathrm{trace}(g|_{\mathrm{ker}P}) - \mathrm{trace}(g|_{\mathrm{coker}P})$. 
This function $L$ is called the 
    \emph{Lefschetz number} of $P$ with respect to the $G$-action. We say the operator $P$ is 
    \emph{rigid for the $G$-action} if its Lefschetz number $L$ is constant as a function on $G$.

One application of rigidity is as follows: 
if we can obtain an expression for the Lefschetz number $L$ of $P$
which does not depend on the choice of $G$-action on $M$ and $W$,
and if we have proved $P$ is rigid in general for such $G$-actions,
then showing the non-constancy of $L$ for a particular 
pair $(W \rightarrow M,P)$
implies that $(W \rightarrow M,P)$ in fact does not admit an equivariant $G$-action. 

In \cite{WittenDiracLoop}, Witten conjectured the rigidity of the twisted Dirac operator $\mathscr{D}$ on a certain element of $K(M)[[q^{1/2}]]$ for $M$ a spin manifold. This $\mathscr{D}$ is viewed as 
an analog of 
the classical Dirac operator
on the (free) loop space of $M$, and the rigidity of this $\mathscr{D}$
is the loop-space analog of the classically known rigidity of the Dirac operator on a finite-dimensional spin manifold. 

Witten's rigidity conjectures were first proved in the work of 
    Hirzebruch \cite{H}, Taubes \cite{T}, Bott and Taubes \cite{BT}, and Krichever \cite{Kr}. Next, in \cite{KL-Mod} and \cite{KL-EGTF}, Liu introduced new methods that simplified and generalized the proof. His approach uses the modularity properties of the Jacobi theta functions to circumvent the technical arguments used in \cite{H} and \cite{BT}.
Subsequently, Dessai \cite{extra1}
extended the work to the spin${}^c$ case,
Liu and Ma \cite{extra2,extra4} explored extensions to families of operators.
Further studies on the odd-dimensional case
were conducted by
Liu and Wang \cite{extra6},
as well as Han and Yu \cite{extra7}.

Let $M$ be a closed even-dimensional spin manifold with an $S^1$-action and $V$ a real even-rank  $S^1$-equivariant spin vector bundle
on $M$. 
Let $\Phi_\lambda$ (actually an element of $K(M)[[q^{1/2}]]$, where $q^{1/2} = e^{\pi \im \tau}$ where $\tau \in \mathbb{H}$ is a complex variable in the upper half-plane), $\lambda=1,2,3$,
be one of the bundles 
considered in the paper \cite{GuanWangII} of Guan and Wang.
The bundle $\Phi_\lambda$ is constructed
from $TM$ and $V$ using standard operations such as direct sums, tensor products,
and symmetric and exterior powers.
Let $\mathscr{D}_\lambda$ be the twisted Dirac operator on (the spinor vector bundle of) $\Phi_\lambda$.
The goal of this paper is to prove 
the following theorem:

\begin{thm}
\label{main-rigidity-thm}


The operators $\mathscr{D}_\lambda$, $\lambda=1,2,3$, 
are rigid if
the first $S^1$-equivariant Pontryagin class of $V$
vanishes.
%
\end{thm}

The proof of Theorem \ref{main-rigidity-thm}
is given in 
    Sections \ref{section:calculatingL} and \ref{section:proofmethod1},
    and has the following outline.
Following the method of \cite{KL-EGTF}, we find that
if we assume the first $S^1$-equivariant Pontryagin class of $V$ vanishes, 
under the modular action, the Lefschetz numbers 
$L_\lambda$ for the $\mathscr{D}_\lambda$ actually transform into each other, up to possible multiplication by nowhere-vanishing holomorphic functions. Following \cite{KL-EGTF}, this ``modularity'' property of the $L_\lambda$'s, together with an analysis of their periodicity properties and possible pole locations, is used to prove the constancy of the $L_\lambda$'s, hence the rigidity of the $\mathscr{D}_\lambda$'s.



\Section{Definitions and conventions}

\begin{conv} 

In this paper, unless otherwise stated, all manifolds are smooth, compact, connected, and without boundary, all vector bundles are smooth and of finite rank and are real except when explicitly declared to be complex, and all actions by Lie groups on manifolds are smooth.

\end{conv}

\begin{notn} ${}$

\begin{itemize}
    \item We use $\im$ to denote the imaginary unit, as opposed to $i$ which is used for just indices and other variables. We denote the $i$th Chern class as $\mathbf{c}_i$, as opposed to $c_i$ which will be used for Chern roots.
    \item We denote by $\xi_{G}$ the ``$G$-equivariant version'' of a given characteristic class $\xi$.
    \item 
    Let $\mathbf{H} := \{ \tau \in \mathbb{C} \,|\, \mathrm{Im}(\tau)>0\}$.
    We let $\tau \in \mathbf{H}$, and we define $q := e^{2\pi \im \tau}$
    and $q^{r} := e^{ 2\pi r \im  \tau }$ for any positive real number $r$.
    \item A formal difference of vector bundles $V$ and $W$ will be denoted $V-W$.
\end{itemize}

\end{notn}

\begin{defin}

Let $\pi : V \rightarrow M$ be a vector bundle. A \emph{$G$-equivariant action on $V \rightarrow M$} consists of:
\begin{itemize}
    \item a $G$-action $\tilde{\rho} : G \times V \rar V$ that linearly maps fibers to fibers, and
    \item the quotient $G$-action $\rho : G\times M \rar M$,
    which satisfies $\pi \circ \tilde{\rho} = \rho \circ (\mathrm{id}_G \times \pi) $.
\end{itemize}

\end{defin}

Let $V \rightarrow M$ be a $G$-equivariant vector bundle.
Let $\Gamma(V)$ denote the space of smooth sections of $V$ and let $\mathscr{D} : \Gamma(V) \rightarrow \Gamma(V)$ be a linear differential operator
which is $G$-equivariant.
Then the kernel and cokernel of $\mathscr{D}$ are representations of $G$.
Suppose further that $\mathscr{D}$ is also Fredholm (e.g.\ an elliptic operator), i.e., the kernel and cokernel of $\mathscr{D}$ are finite-dimensional. 

\begin{defin}
    \label{def:Lef-num}

In the above situation, the \emph{Lefschetz number} of $\mathscr{D}$
is the function
$ L = L_{\mathscr{D}} : G \rar \mathbb{F}$ given by 
$$
L_{\mathscr{D}}(g) := \mathrm{trace}(g|\mathrm{ker}(\mathscr{D})) - \mathrm{trace}(g|\mathrm{coker}(\mathscr{D}))
,
$$
where we use the convention that the codomain $\mathbb{F}$ of $L_{\mathscr{D}}$ is $\bb{R}$ or $\bb{C}$ according to whether $V$ is real or complex.  

\end{defin}

\begin{rem}
    Later, 
    instead of a single vector bundle $V$ on $M$ with a differential operator $\mathscr{D}$, we will consider a formal power series 
    $ \sum_{n=0}^{\infty} V_{(n)} q^n $ where $q$ is a formal variable and 
    each $V_{(n)}$ is itself an $S^1$-equivariant vector bundle
    with an $S^1$-equivariant linear Fredholm differential operator $\mathscr{D}_{(n)}$.
    In this case, the ``Lefschetz number'' $L$ of the collection $(\mathscr{D}_{(n)})_{n=0,1,2,\ldots}$
    is the function $L : G \rar \mathbb{F}[[q]]$
    given by $L(g) = \sum_{n=0}^{\infty} L_n(g)q^n$
    where $L_n$ is the Lefschetz number of $\mathscr{D}_{(n)}$ for each $n\geq0$.
    Again $\mathbb{F}$ is $\bb{R}$ or $\bb{C}$ according to whether the $V_n$'s are real or complex.  
\end{rem}

\begin{defin}

The operator $\mathscr{D}$ is $G$-\textit{rigid}
if the Lefschetz number $L$ of $\mathscr{D}$ is constant on $G$, which is the case 
if and only if for every $g \in G$
we have 
$$
L(g) = L(1_G) = \mathrm{dim}\,\mathrm{ker}(\mathscr{D}) - \mathrm{dim}\,\mathrm{coker}(\mathscr{D}).
$$
We note the RHS is just the ordinary Fredholm index of $\mathscr{D}$.



\end{defin}

\begin{rem} [Lefschetz number is well-defined on $K$-group]

Given $G$-equivariant vector bundles $A,B,A',B'$ on $M$ with respective $G$-equivariant differential operators $\mathscr{D}_A,\mathscr{D}_B,\mathscr{D}_{A'},\mathscr{D}_{B'}$, if there is an isomorphism 
$A \oplus B' \xrightarrow{\simeq} A' \oplus B$ of $G$-equivariant bundles that also identifies $\mathscr{D}_A \oplus \mathscr{D}_{B'}$ with $\mathscr{D}_{A'} \oplus \mathscr{D}_B$, 
then we have $\mathscr{L}_A + \mathscr{L}_{B'} = \mathscr{L}_{A'} + \mathscr{L}_B $ for the corresponding Lefschetz numbers,
so that $\mathscr{L}_A - \mathscr{L}_B = \mathscr{L}_{A'} - \mathscr{L}_{B'}$.

Therefore, the Lefschetz number of a formal difference 
$A - B$
of ``concrete'' $G$-equivariant vector bundles $A,B$ both equipped with $G$-equivariant linear differential operators, 
is viewed as just the difference of their Lefschetz numbers, $\mathscr{L}_A - \mathscr{L}_B$.

\end{rem}

\Section{Reduction to the case of $G=S^1$}

Our aim in this paper is to prove the rigidity of twisted Dirac operators $\mathscr{D}$ on certain bundles,
when equivariant for an action by any Lie group $G$ that is compact and connected.

The goal of this section is to reduce to the case of $G = S^1$.
First we need the following lemma:

\begin{lem}
    \label{lem:Lef-is-continuous}

The Lefschetz number $L : G \rightarrow \mathbb{F}$
of a $G$-equivariant Fredholm operator $\scr{D}$ on a $G$-equivariant vector bundle $V \rightarrow M$
is continuous as a function on $G$. 

\begin{proof}

Since the trace is continuous on the space of endomorphisms of a finite-dimensional vector space (over $\mathbb{F}$), 
it suffices to show 
    the action homomorphism $G \rightarrow \mathrm{GL}_\bb{F}(A)$
is continuous 
    when $A$ is either $\mathrm{ker}(\scr{D})$ or $\mathrm{coker}(\scr{D})$.

Let $C(M)$ denote the space of all continuous $\mathbb{F}$-valued functions on $M$. We give $C(M)$ the compact-open topology.

In general we will show the action homomorphism $G \rightarrow \mathrm{GL}_\bb{F}(A)$ is continuous for $A$ any $\mathbb{F}$-linear finite-dimensional subspace of $C(M)$ which is stable under the $G$-action. 

If $G$ acts continuously on $M$, then also $G$ acts continuously on $C(M)$ with the compact-open topology. 
So, if $\mathscr{F}$ denotes the space of all continuous $\mathbb{F}$-linear automorphisms of $C(M)$, with $\mathscr{F}$ itself having the compact-open topology
(inherited from the space of all continuous maps $C(M) \rightarrow C(M)$),
we have a continuous action homomorphism $G \rightarrow \mathscr{F}$.

Let $A \subset C(M)$ be a given $G$-stable finite-dimensional linear subspace.
Let $\mathscr{F}_A \subset \mathscr{F}$ be the set
of all $f \in \mathscr{F}$ such that $f$ maps $A$ linearly onto itself.
Now the action homomorphism $G \rightarrow \mathscr{F}_A$ with codomain restricted to $\mathscr{F}_A$
is still continuous, with $\mathscr{F}_A$ having the subspace topology from $\mathscr{F}$.

Next, $\mathrm{GL}_\bb{F}(A)$ 
can be identified topologically with the quotient of $\mathscr{F}_A$
by the relation considering two functions $h,k \in \mathscr{F}_A$
equivalent if $h,k$ agree on $A \subset C(M)$.

Now the composition 
    $G \rightarrow \mathscr{F}_A \xrightarrow{\text{quotient}} \mathrm{GL}_\bb{F}(A)$ 
is continuous and agrees with the action homomorphism $G \rightarrow \mathrm{GL}_{\mathbb{F}}(A)$. Thus $G \rightarrow \mathrm{GL}_{\mathbb{F}}(A)$ is continuous.
\end{proof}

\end{lem}

We also recall the following facts about
Lie groups:

\begin{fact}
\label{lem:U-dense-Abel}

    Let $T = \mathbb{R}^n/\mathbb{Z}^n$.
    Let $U \subset T$
    be the set of elements $h \in T$
    such that there exists a closed Lie subgroup $H \subset T$
    with $h \in H$
    and such that $H \cong S^1$ as a Lie group.
    Then $U$ is dense in $T$.



    
\end{fact}

\begin{lem}
\label{lem:U-dense-general-cptconn}

Let $G$ be any compact connected Lie group.
    Let $U \subset G$
    be the set of elements $h \in G$
    such that there exists a closed Lie subgroup $H \subset G$
    with $h \in H$
    and such that $H \cong S^1$ as a Lie group.
    Then $U$ is dense in $G$.

\begin{proof}
Let $O \subset G$ be any nonempty open subset.
We aim to show $O \cap U \neq \emptyset$.

Pick any $g \in O$.
We know (see for example Corollary 4.46 of \cite{Knapp}) there 
exists a maximal torus $T \subset G$
such that $g \in T$. 
By Fact \ref{lem:U-dense-Abel},
there exists an $h \in O \cap T$ 
such that $h$ is contained in some closed Lie subgroup $H \subset T$ with $H$ isomorphic to $S^1$. But now $h \in U$, by definition of $U$.

Hence $h \in O \cap U$, proving $U \cap O \neq \emptyset$.
\end{proof}

\end{lem}

Now we prove the main result of this section:

\begin{lem}
    If an operator $\mathscr{D}$ is rigid with respect to any $S^1$-action,
    then it is rigid with respect to the action of any compact connected Lie group.

\begin{proof}

Suppose $\mathscr{D}$ is an operator on some pair $(W\rightarrow M,P)$ where $P$ is a differential operator on a vector bundle $W\rightarrow M$, and assume $\mathscr{D}$ is rigid with respect to all $S^1$-actions. 
Let $G$ be a compact connected Lie group with a given equivariant action on $(W\rightarrow M,P)$.
Let $L = L_\mathscr{D}$ be the Lefschetz number of $\mathscr{D}$ with respect to the $G$-action.
We aim to show $\mathscr{D}$ is rigid, i.e.\ that $L$ is constant,
or equivalently that $L(g) = L(1_G)$ for any $g \in G$.

We know $L$ is continuous on $G$ by Lemma \ref{lem:Lef-is-continuous}.
Therefore, it suffices to show that $L$ takes the value $L(1_G)$
on a dense subset of $G$.

Let $U \subset G$ be the subset consisting of all $g \in G$
such that there exists a Lie subgroup $H \subset G$
such that $H \cong S^1$ and $g \in H$. 
Since $G$ is compact and connected, by Lemma \ref{lem:U-dense-general-cptconn} we know $U$ is dense in $G$.

We now show for any given $g \in U$, $L(g) = L(1_G)$. Given $g \in U$, find $H \subset G$ an isomorphic copy of $S^1$ containing $g$. By hypothesis, $\mathscr{D}$ is rigid for the $H$-action, so it follows that $L(g) = L(1_H) = L(1_G)$, as was to be shown.
\end{proof}

\end{lem}

Therefore, in what follows, we will consider only actions by the circle group $G = S^1$.

\Section{Tools for calculating Lefschetz numbers}

    Let $V$ be an oriented rank $2k$ real $S^1$-equivariant vector bundle on a manifold $M$, where the $S^1$-action on $M$ is trivial.
    Assume $M$ is connected.
    Then there exist $m_1,\ldots,m_k \in \mathbb{Z}$
    such that for every $p \in M$, the fiber $V_p$ of $V$ at $p$
    decomposes as a sum of real 2-planes $V_p =  L_1 \oplus \cdots \oplus L_k $
    where each $L_i$
    is $S^1$-stable with weight $m_i$ for the $S^1$-action.


\begin{defin} \label{def:weights}
    The integers $m_1,\ldots,m_k$ 
        above
    are called the \emph{weights}
    of the $S^1$-action on $V$.
\end{defin}


Let $w \in H^2_{S^1}(M)$ denote the ``universal element''
given by the pullback of $\widetilde{w} \in H^2_{S^1}(\{\text{pt}\})$ along the map $M \rightarrow \{\text{pt}\}$,
where $\widetilde{w} \in H^2_{S^1}(\{\text{pt}\}) \cong H^2(BS^1)$ 
is the ``universal first Chern class.'' 

\begin{lem}
\label{lem:equivariant-chern-roots}

The formal $S^1$-equivariant Chern roots of 
the complexification
$V_{\bb{C}}$ of $V$ are
$$
\left\{ \pm ( \xi_i + m_i w) \right\}_{\pm; i=1,\ldots,k}
$$
with $\left\{ \pm \xi_i \right\}_{\pm; i=1,\ldots,k}$ the formal (non-equivariant) Chern roots
of $V_\bb{C}$.

\end{lem}

The $\hat{A}$-class $\hat{\mathscr{A}}(V_\bb{C})$, the Chern character $\mathrm{ch}(V_\bb{C})$, the Euler class $E(V)$, and the total Pontryagin class $p(V)$ admit (see e.g. \cite{AS3}) 
the following expressions in terms of the formal Chern roots $\pm \xi_1,\ldots,\pm\xi_k$ of $V_\bb{C}$:
\begin{equation} \label{eqn:nonequiv-classes}
\begin{matrix}
    \hat{\mathscr{A}}(V) = \prod_{i=1,\ldots,k} [f(\xi_i)f(-\xi_i) ] = \prod_{i=1,\ldots,k} g(\xi_i) 
    ,&
    \mathrm{ch}(V_{\bb{C}}) 
        = \sum_{i=1,\ldots,k} 
        2\cosh ( \xi_i) ,
        \\
    E(V) = \prod_{i=1,\ldots,k} \xi_i                                               ,&
    p(V) = \prod_{i=1,\ldots,k} ( 1+\xi_i^2 ) ,
\end{matrix}
\end{equation}
where $f(z),g(z)$ are given by the Maclaurin series of the functions 
$$ \begin{matrix}
f(z) = \sqrt{g(z)} ,& g(z) = \frac{z/2}{\mathrm{sinh}(z/2)} .
\end{matrix}
$$

The $S^1$-equivariant versions 
of the above characteristic classes 
are obtained 
by replacing the formal Chern roots $\pm \xi_i$ in the above expressions
by the corresponding formal $S^1$-equivariant Chern roots. Hence: 

\begin{lem}
\label{lem:equivariant-classes}

In the notation of Lemma \ref{lem:equivariant-chern-roots},
the $S^1$-equivariant versions of the characteristic classes 
in Equation \ref{eqn:nonequiv-classes}
are
\begin{equation*}
\begin{matrix}
    \hat{\mathscr{A}}_{S^1}(V)
    = \prod_{i=1,\ldots,k} g(\xi_i + m_i w)
    ,&
    \mathrm{ch}_{S^1}(V_{\bb{C}}) 
        = \sum_{i=1,\ldots,k}    
        2\cosh ( \xi_i + m_i w) ,
        \\
    E_{S^1} (V) = \prod_{i=1,\ldots,k} (\xi_i + m_i w)                                            ,&
    p_{S^1}(V) = \prod_{i=1,\ldots,k} ( 1+(\xi_i + m_i w)^2 ) ,
\end{matrix}
\end{equation*}
where as before, $g(z) = (z/2)/\sinh(z/2)$, and $w \in H^2_{S^1}(M)$ is the universal element.

\end{lem}

Let $\mathscr{D}$ be the Dirac operator on a spin manifold $M$ equipped with an $S^1$-action, now no longer necessarily the trivial action. 
Let $V$ be an $S^1$-equivariant spin vector bundle on $M$.
Let $L : S^1 \rightarrow \mathbb{C}$ be the Lefschetz number for the ``twisted Dirac operator'' $\mathscr{D} \otimes V$ on (the spinor vector bundle of) $V$.

\begin{lem}
    \label{lem:formula-for-Lef}
For $g \in S^1$ sufficiently near the identity we have
$$
L_{\mathscr{D}}(g) = 
\sum_{\alpha \in \mathcal{I}} 
\int_{M_\alpha}
\left. 
\left( \frac {  \hat{\mathscr{A}}_{S^1}(i_\alpha^* TM) \,\mathrm{ch}_{S^1}(i_\alpha^* V_\bb{C} )  }{E_{S^1}(\nu_\alpha)} \right)
\right|_{w = 2\pi \im t_g} 
$$
where:
\begin{itemize}
    \item $\{M_\alpha\}_{\alpha \in \mathcal{I}}$ lists the connected components of the fixed locus (set of fixed points) of the $S^1$-action, and for each $\alpha \in \mathcal{I}$, 
    $\nu_\alpha$ is the normal bundle of $M_\alpha$ in $M$, 
    and $i_\alpha : M_\alpha\hookrightarrow M$ is the inclusion;

\item $t_g \in \bb{R}$ is sufficiently near 0 
such that $e^{2\pi \im t_g} = g$, 
$g \in S^1$;

\item $w$ is the canonical universal element of $H^2_{S^1}(M)$;
and

\item the notation $(\cdots)|_{w=2\pi \im t_g}$ means that we take the expression $(\cdots)$ and evaluate $w$ as $2\pi \im t_g$.

\end{itemize}

\begin{proof}
    This is a consequence of 
the Atiyah-Singer-Lefschetz formula 
    (see e.g.\ Section 2 of \cite{AS3})
and Bott localization.
\end{proof}

\end{lem}

We note that the fixed-point set of the $S^1$-action
on $M$ is a smooth submanifold of $M$,
since a smooth proper action by a Lie group $G$ on a manifold $M$ always has fixed-point set $M^G$ a smooth submanifold of $G$, and $S^1$ is compact and acting continuously on a Hausdorff space $M$
so the action is indeed proper.

\begin{rem}
    \label{rem:inf-prod-exprs-for-theta-funcs}


We will use the following infinite product expressions for the Jacobi theta functions,
which can be found 
in Section 2.2 of \cite{GuanWangII}:
$$
\theta(v,\tau) = c \cdot q^{1/8} \cdot 2\sin(\pi v)
\cdot \prod_{r=1}^{\infty}
 \left[ (1-q^r e^{2\pi \im v}) (1-q^r e^{-2\pi \im v})  \right] 
 , $$$$
\theta_1(v,\tau) = c \cdot q^{1/8} \cdot 2\cos(\pi v)
\cdot \prod_{r=1}^{\infty}
 \left[ (1 + q^r e^{2\pi \im v}) (1 + q^r e^{-2\pi \im v})  \right] 
 , $$$$
\theta_2(v,\tau) = c 
\cdot \prod_{r=1}^{\infty}
 \left[ (1 - q^{r-1/2} e^{2\pi \im v}) (1 - q^{r-1/2} e^{-2\pi \im v})  \right] 
 , $$$$
\theta_3(v,\tau) = c 
\cdot \prod_{r=1}^{\infty}
 \left[ (1 + q^{r-1/2} e^{2\pi \im v}) (1 + q^{r-1/2} e^{-2\pi \im v})  \right] 
$$
where $q^{1/8} = e^{\pi \im \tau / 4}$, and
$$
c = \prod_{r=1}^{\infty} (1-q^r)  .
$$

We will also use the function $\theta'(0,\tau) := \left. \frac{\partial}{\partial v}\right|_{v=0} \theta(v,\tau)$.
Using the Jacobi identity (see Theorem 5 in Section 5 of \cite{Chand}), which is
$$
\theta'(0,\tau) = \pi \theta_1(0,\tau)\theta_2(0,\tau)\theta_3(0,\tau),
$$
we also obtain the following expression for $\theta'(0,\tau)$:
$$
\theta'(0,\tau) = 2\pi q^{1/8} \prod_{r=1}^{\infty} \left[ (1-q^r)^3 \right]   .
$$


\end{rem}

Let $W$ 
    be a possibly virtual, arbitrary complex vector bundle of complex rank $k$, on a base space $M$, with formal Chern roots $ \xi_1,\ldots, \xi_k$.
    
\begin{defin}

Let $K(M)[[t]]$ be the ring of formal power series in a variable $t$ with $K(M)$ coefficients. 
We make the standard definitions $S_t W := \oplus_{r=0}^{\infty} t^r S^r W $ and $\Lambda_t W := \oplus_{r=0}^{\infty} t^r \Lambda^r W $,
where $S^r,\Lambda^r$
denote the $r$th symmetric and exterior powers respectively.
These may be viewed as elements of $K(M)[[t]]$.

\end{defin}

We recall the following well-known formulas for Chern characters (see e.g.\ \cite{GuanWangII}). 
Analogous formulas hold for the below equivariant Chern characters 
when $W$ is a complexification of a real even-rank bundle with fiberwise $S^1$-action 
by replacing the non-equivariant Chern roots by the equivariant ones.

\begin{lem} 
    \label{lem:ch-exprs}

If $\xi_1,\ldots,\xi_k$ are the formal Chern roots of $W$, then:

\begin{itemize}
    
    \item $\mathrm{ch}(S_t W) = \prod_{i=1,\ldots,k} \left( \frac{1}{1 - te^{\xi_i}} \right)  $ ; 
$\mathrm{ch}(\Lambda_t W) = \prod_{i=1,\ldots,k} (1 + te^{\xi_i})  $;

    \item $\mathrm{ch}(S_t(-W)) = \prod_{i=1,\ldots,k} ( 1 - te^{\xi_i} ) $ ; $\mathrm{ch}(\Lambda_t(-W)) = \prod_{i=1,\ldots,k} \left( \frac{1}{1 + te^{\xi_i}} \right) $;

    \item if $\epsilon_N$ denotes the trivial rank-$N$ complex vector bundle, then $\mathrm{ch}(S_t \epsilon_N) = \frac{1}{(1-t)^N}$ and $\mathrm{ch}(\Lambda_t \epsilon_N) = (1+t)^N$.

    
\end{itemize}

We denote by $\widetilde{W}$ the formal difference $W - \mathrm{dim}\,W$, where $\mathrm{dim}\, W$ represents a trivial complex bundle of the same complex rank as $W$.

Next, let us assume $W$ is the complexification of a real rank $k=2\ell$ bundle $V$. 
Then the Chern roots of $W$ can be taken to be $\pm 2\pi \im c_i$ for each sign $\pm$ and $i=1,\ldots,\ell$.
Here the extra $2\pi \im $ factors were added because they will be useful later.

Recall that $\tau\in \mathbf{H}$, $q = e^{2\pi \im \tau}$, and $q^{1/2} = e^{\pi \im \tau}$.
Then by the multiplicativity of the Chern character, 
the formulas in Lemma \ref{lem:ch-exprs},
and the infinite product expressions for the Jacobi theta functions
in Remark \ref{rem:inf-prod-exprs-for-theta-funcs},
the following hold:

\end{lem}

\begin{lem}
    \label{lem:ch-inf-tensor-prods}

$$
\begin{matrix}
    \mathrm{ch}\left( \bigotimes_{r=1}^{\infty} S_{q^r}    
        (\widetilde{W}) \right) 
    = \prod_{i=1,\ldots,\ell} \left[  \frac{\sin(\pi c_i)}{\pi} \frac{\theta'(0,\tau)}{\theta(c_i,\tau)}  \right]
    ,&
    \mathrm{ch}\left( \bigotimes_{r=1}^{\infty} \Lambda_{q^r}  
        (\widetilde{W}) \right) 
    = \prod_{i=1,\ldots,\ell} \left[  \frac{1}{\cos(\pi c_i)} \frac{\theta_1(c_i,\tau)}{\theta_1(0,\tau)} \right]
    ,\\
    \mathrm{ch}\left( \bigotimes_{r=1}^{\infty} S_{q^{r-1/2}} 
        (\widetilde{W}) \right) 
    = \prod_{i=1,\ldots,\ell} \left[  \frac{\theta_2(0,\tau)}{\theta_2(c_i,\tau)}  \right]
    ,&
    \mathrm{ch}\left( \bigotimes_{r=1}^{\infty} \Lambda_{ - q^{r-1/2}} 
        (\widetilde{W}) \right) 
    = \prod_{i=1,\ldots,\ell} \left[  \frac{\theta_2(c_i,\tau)}{\theta_2(0,\tau)}  \right]
    ,\\
    \mathrm{ch}\left( \bigotimes_{r=1}^{\infty} S_{ - q^{r-1/2}} 
        (\widetilde{W}) \right) 
    = \prod_{i=1,\ldots,\ell} \left[  \frac{\theta_3(0,\tau)}{\theta_3(c_i,\tau)}  \right]
    ,&
    \mathrm{ch}\left( \bigotimes_{r=1}^{\infty} \Lambda_{q^{r-1/2}} 
        (\widetilde{W}) \right) 
    = \prod_{i=1,\ldots,\ell} \left[  \frac{\theta_3(c_i,\tau)}{\theta_3(0,\tau)}  \right]
    .
\end{matrix}
$$

To get the ``$S^1$-equivariant versions'' of the above formulas (which hold when $V$ has an $S^1$-action that preserves each fiber),
we simply replace each $c_i$ by $c_i + m_i w/(2\pi\im)$,
where $w$ is the ``universal element'' of $H^2_{S^1}$ 
as described above Lemma \ref{lem:equivariant-chern-roots}
and $m_1,\ldots,m_\ell$ are the exponents of the $S^1$-action on $V_0$.


\end{lem}

\Section{Calculation of the Lefschetz numbers $\mathscr{L}$ of the twisted Dirac operators $\mathscr{D}$ on the bundles $\Phi$ considered by Guan and Wang}
    \label{section:calculatingL}

Let $M$ be a manifold. Following Guan and Wang \cite{GuanWangII}, 
    if $V$ is a spin vector bundle on $M$,
then we denote by $\Delta(V)$ the spinor bundle corresponding to $V$. 
    (Specifically, $\Delta(V)$ is the vector bundle 
    associated to the principal $\mathrm{Spin}(\mathrm{rk}V)$ bundle of $V$ via the spinor representation of $\mathrm{Spin}(\mathrm{rk}V)$.)
If $M$ itself is a spin manifold, then $TM$ is spin and we 
simply denote $\Delta(TM)$ by $\Delta(M)$. 
Furthermore, if $V$ has real even rank $2\ell$
and $V_\bb{C}$ has formal Chern roots $\{ \pm 2\pi \im c_i \}_{\pm;i=1,\ldots,\ell }$,
then
\begin{equation}
    \label{eqn:spinor-bundle-ch}
\mathrm{ch}(\Delta(V)) = \prod_{i=1,\ldots,\ell} (e^{\pi \im c_i} + e^{-\pi \im c_i})
= 2^\ell \prod_{i=1,\ldots,\ell} \cos(\pi c_i) .
\end{equation}

For example, 
by Lemma \ref{lem:ch-inf-tensor-prods}
and Equation \ref{eqn:spinor-bundle-ch},
$$
    \mathrm{ch}\left( \Delta(V) \otimes \bigotimes_{r=1}^{\infty} \Lambda_{q^r} ( \widetilde{V_\bb{C}} ) \right)
= 2^\ell \prod_{i=1,\ldots,\ell} \left[  \frac{\theta_1(c_i,\tau)}{\theta_1(0,\tau)}  \right] .
$$

\vspace{1em}

We introduce the following bundles
(or more precisely power series in $q^{1/2}$ with coefficients in $K(M)$)
 considered by Guan and Wang (\cite{GuanWangII}). 


For $W$ a given real vector bundle on $M$, we define the following elements of $K(M)[[q^{1/2}]]$ :
\begin{equation}
    \label{def:thetas}
\begin{matrix}
    \Theta_1  (W) := \left( \bigotimes_{r=1}^{\infty} S_{q^r} (\widetilde{W}) \right) \otimes \left( \bigotimes_{s=1}^{\infty} \Lambda_{q^s} (\widetilde{W}) \right)    ,\\
     \Theta_2  (W) := \left( \bigotimes_{r=1}^{\infty} S_{q^r} (\widetilde{W}) \right) \otimes \left( \bigotimes_{s=1}^{\infty} \Lambda_{-q^{s-1/2}} (\widetilde{W}) \right)    ,\\
     \Theta_3 (W) := \left( \bigotimes_{r=1}^{\infty} S_{q^r} (\widetilde{W}) \right) \otimes \left( \bigotimes_{s=1}^{\infty} \Lambda_{q^{s-1/2}} (\widetilde{W}) \right)   ,
\end{matrix}
\end{equation}
where as in Lemma \ref{lem:ch-inf-tensor-prods}
we are using the notation 
$
\widetilde{U} := U - \mathrm{dim}\,U
$ for any complex/real vector bundle $U$;
here $\mathrm{dim}\,U$ is a trivial complex/real bundle of same complex/real rank as $U$.
(Also, in the $G$-equivariant situation, 
the action of $G$ on the total space of a trivial vector bundle $\epsilon$ is taken to be the one that lies over the given action on the base space $M$
    and sends each fiber to its image via just an ``identity'' map, with respect to some global trivialization of $\epsilon$.)

If $\pm 2\pi \im c_i$ for $i=1,\ldots,k$ is a list of formal Chern roots of $V_\bb{C}$ for $V$ a real rank $2k$ bundle, then 
    by Lemma \ref{lem:ch-inf-tensor-prods} and Equation \ref{def:thetas},
\begin{equation}
\begin{matrix}
    \mathrm{ch}((\Theta_1(V))_\bb{C}) = \frac{1}{\pi^\ell} \prod_{i=1,\ldots,\ell} \left[ \tan(\pi c_i) 
        \frac{ \theta'(0,\tau) \theta_1(c_i,\tau)} { \theta(c_i,\tau) \theta_1(0,\tau) } \right]
    ,\\
   \mathrm{ch}((\Theta_2(V))_\bb{C}) = \frac{1}{\pi^\ell} \prod_{i=1,\ldots,\ell} \left[ \sin(\pi c_i) 
        \frac{ \theta'(0,\tau) \theta_2(c_i,\tau)} { \theta(c_i,\tau) \theta_2(0,\tau) } \right]
    ,\\
   \mathrm{ch}((\Theta_3(V))_\bb{C}) = \frac{1}{\pi^\ell} \prod_{i=1,\ldots,\ell} \left[ \sin(\pi c_i) 
        \frac{ \theta'(0,\tau) \theta_3(c_i,\tau)} { \theta(c_i,\tau) \theta_3(0,\tau) } \right]
    .
\end{matrix}
\end{equation}

Next, let $M$ be a spin manifold of dimension $2d$, and let $V$ be a real rank $2\ell$ spin vector bundle on $M$.
Following \cite{GuanWangII} we further define the following elements of $K(M)[[q^{1/2}]]$:
\begin{equation}
    \label{def:Phis}
\begin{matrix}
    \Phi_0(M) := \big(\Delta(M) \otimes (\Theta_1(TM))_\bb{C} \big) \oplus \left[  (\Theta_2(TM))^{\oplus 2^d} \oplus  (\Theta_3(TM))^{\oplus 2^d} \right]_\bb{C} 
        ,\\
    \Phi_1(M,V) := \Phi_0(M) \otimes \Delta(V) \otimes \bigotimes_{r=1}^{\infty} \Lambda_{q^r} (\widetilde{V_\bb{C}})  
        ,\\
    \Phi_2(M,V) := \Phi_0(M) \otimes \bigotimes_{r=1}^{\infty} \Lambda_{-q^{r-1/2}} (\widetilde{V_\bb{C}})  
        ,\\
    \Phi_3(M,V) := \Phi_0(M) \otimes \bigotimes_{r=1}^{\infty} \Lambda_{q^{r-1/2}} (\widetilde{V_\bb{C}})  
        .
\end{matrix}
\end{equation}

Applying Lemma \ref{lem:ch-inf-tensor-prods}
and Equation \ref{eqn:spinor-bundle-ch}
    to Equation \ref{def:Phis},
    we see
    that if $\pm 2\pi \im b_j$ for $j=1,\ldots,\ell$ is the list of Chern roots for $V_\bb{C}$,
    and $\pm 2\pi \im c_i$ for $i=1,\ldots,d$ the list of Chern roots for $TM_\bb{C}$,
    then
\begin{equation}
    \label{eqn:Phis-chs}
\begin{matrix}
\mathrm{ch}(\Phi_0) = \left( \tfrac{2}{\pi} \right)^d \left( \prod_{i=1,\ldots,d} \sin(\pi c_i) \right)
\left( \sum_{\mu=1}^{3}
\prod_{i=1,\ldots,d} \frac{ \theta'(0,\tau) \theta_\mu(c_i,\tau)} { \theta(c_i,\tau) \theta_\mu(0,\tau) } 
\right) ,\\
\mathrm{ch}(\Phi_1) = 
2^\ell
(\mathrm{ch}\Phi_0) 
\left( \prod_{j=1}^{\ell}   \frac{\theta_1(b_j,\tau)}{\theta_1(0,\tau)} \right)
,\\
\mathrm{ch}(\Phi_2) = (\mathrm{ch}\Phi_0)  
\left( \prod_{j=1}^{\ell}   \frac{\theta_2(b_j,\tau)}{\theta_2(0,\tau)} \right)
,\\
\mathrm{ch}(\Phi_3) = (\mathrm{ch}\Phi_0) 
\left( \prod_{j=1}^{\ell}   \frac{\theta_3(b_j,\tau)}{\theta_3(0,\tau)} \right)
.
\end{matrix}
\end{equation}

Let $\mathscr{D}$ be the Dirac operator on $M$.
Now for each $\lambda =1,2,3$ put $\mathscr{D}_\lambda := \mathscr{D} \otimes \Phi_\lambda$, the twisted Dirac operator on $\Phi_\lambda$, as a differential operator
on the spinor vector bundle of $\Phi_\lambda$.
Let us assume $V \rightarrow M$ is an $S^1$-equivariant bundle,
for which $\mathscr{D} \otimes V$ is also equivariant (on the spinor bundle of $V$).
This results in an action on $\Phi_\lambda$, for which $\mathscr{D}_\lambda = \mathscr{D}\otimes \Phi_\lambda$ is equivariant (on the spinor bundle of $\Phi_\lambda$).

Applying Lemmas \ref{lem:equivariant-classes} and \ref{lem:formula-for-Lef}
and
 Equation \ref{eqn:Phis-chs},
for $g \in S^1$ sufficiently near the identity, 
for $\lambda =1,2,3$ we have the following expression
for the Lefschetz number $\mathscr{L}_{\mathscr{D}_\lambda}$
of $\mathscr{D}_\lambda$:
\begin{lem} 
    \label{lem:lef-expression}
\begin{multline*} 
       \mathscr{L}_{\mathscr{D}_\lambda}( g ; \tau) =  ( \kappa_\lambda ) \left( \tfrac{2}{\pi} \right)^d  
       \sum_{\alpha\in\mathcal{I}} 
            \int_{M_\alpha}   
\left[
    \left( \prod_{i=1}^{d_\alpha} c_{\alpha,i}
\right)
\right. \\ \left. \times
    \left(
\sum_{\mu=1}^{3}\prod_{i=1}^{d} \frac{ \theta'(0,\tau) \theta_\mu(c_{\alpha,i} + n_{\alpha,i}t ,\tau)} { \theta(c_{\alpha,i} + n_{\alpha,i}t,\tau) \theta_\mu(0,\tau) } 
\right) 
    \left( \prod_{j=1}^{\ell}  
 \frac{\theta_\lambda (b_{\alpha,j} + m_{\alpha,j}t,\tau)}{\theta_\lambda(0,\tau)}
\right)
\right] \;,
\end{multline*}
where:
\begin{itemize}
\item $t \in \bb{R}$ is sufficiently close to 0 and is such that $g = e^{2\pi \im t} \in S^1$.
\item $2d = \mathrm{dim}M$ and $2\ell = \mathrm{rk}V$.
\item 
$\kappa_\lambda=2^{\ell} $ 
if $\lambda=1$, and $\kappa_\lambda=1$ if $\lambda=2,3$.
\item $\{M_\alpha\}_{\alpha\in\mathcal{I}}$ lists the connected components of the fixed locus (set of fixed points) of the $S^1$-action on $M$,
and $d_\alpha = \mathrm{dim}\,M_\alpha$.
\item For each connected component $M_\alpha$ of the fixed locus, the following hold.
\begin{itemize}

\item The weights (see Definition \ref{def:weights}) of $i_\alpha^*TM$
are $n_{\alpha,i}$ for $i=1,\ldots,d$,
with $n_{\alpha,i}=0$ for $1\leq i\leq d_\alpha$
being the weights of $TM_\alpha$, which are just 0.

\item 
The (non-equivariant)
formal Chern roots of $(i_\alpha^*TM)_\bb{C}$
are 
$ \pm 2\pi \im c_{\alpha,i} $ for $i=1,\ldots,d$,
with $ \pm 2\pi \im c_{\alpha,i} $ for $i=1,\ldots,d_\alpha$
being the formal Chern roots of $(TM_\alpha)_\bb{C}$.

\item The weights (see Definition \ref{def:weights}) of $i_\alpha^*V$
are $m_{\alpha,j}$ for $j=1,\ldots,\ell$.

\item 
The (non-equivariant)
formal Chern roots of $(i_\alpha^*V)_\bb{C}$
are
$ \pm 2\pi \im b_{\alpha,j} $ for $j=1,\ldots,\ell$.

\end{itemize}
\end{itemize}

\begin{proof}
    The calculation follows immediately from Lemmas
\ref{lem:equivariant-classes} and \ref{lem:formula-for-Lef}
and Equation \ref{eqn:Phis-chs}.
\end{proof}

\end{lem}

In what follows, we shall view the above expression as defining a function 
$
\mathscr{L}_\lambda(t,\tau) \,,
$
for each $\lambda=1,2,3$.
The expression for $\mathscr{L}_\lambda$ involves rational expressions in the Jacobi theta functions, which are holomorphic everywhere in $\mathbb{C} \times \mathbf{H}$, and so we can extend the domain of $\mathscr{L}_\lambda$ to every complex number $(t,\tau) \in \mathbb{C} \times \mathbf{H}$ where the denominators are nonzero, and $\mathscr{L}_\lambda(t,\tau)$ is then meromorphic on $\mathbb{C} \times \mathbf{H}$.

\vspace{1em}

To prove the rigidity of $\mathscr{D}_\lambda$ on $\Phi_\lambda$, $\lambda = 1,2,3$, it now suffices to show $\mathscr{L}_\lambda(t,\tau)$ is
constant in $t$ (for fixed $\tau \in \mathbf{H}$), at least for real values of $t$.
In the next section (Section \ref{section:proofmethod1}), we do this
``directly'' via a method following that of \cite{KL-EGTF}.

\Section{Analysis of the Lefschetz numbers $\mathscr{L}$ using Liu's techniques involving Jacobi theta function properties}
    \label{section:proofmethod1}

Throughout Section \ref{section:proofmethod1} 
we use the notation $\mathscr{L}_\lambda(t,\tau)$, $\lambda=1,2,3$,
for the meromorphic function (on $\mathbb{C} \times \mathbf{H}$) whose expression 
we found in Lemma \ref{lem:lef-expression}.

\SubSection{Double periodicity of $\mathscr{L}$}  

We will reference the following properties of theta-functions, which can be found in e.g.\ Chapter V of \cite{Chand}:
\begin{equation}
\label{eqn:theta-shifts}
\begin{matrix}
    \theta(v+1,\tau) =  - \theta(v,\tau)  ,&   \theta(v+\tau,\tau) =  - q^{-1/2} e^{-2\pi \im v} \theta(v,\tau) ,\\
    \theta_1(v+1,\tau) =  - \theta_1(v,\tau)  ,&   \theta_1 (v+\tau,\tau) =  +q^{-1/2} e^{-2\pi \im v} \theta_1(v,\tau) ,\\
    \theta_2 (v+1,\tau) =  +\theta_2(v,\tau)  ,&   \theta_2(v+\tau,\tau) =  - q^{-1/2} e^{-2\pi \im v} \theta_2(v,\tau) \\
    \theta_3 (v+1,\tau) =  +\theta_3(v,\tau)  ,&   \theta_3(v+\tau,\tau) =  +q^{-1/2} e^{-2\pi \im v} \theta_3(v,\tau)  .
\end{matrix} 
\end{equation}
(We note the convention in this paper is $q = e^{2\pi \im \tau}$, while \cite{Chand} uses $q = e^{ \pi \im \tau}$.)

By applying the formulas in the 2nd column repeatedly, we obtain the following formulas which are valid for any integer $k$:
\begin{equation}
\label{eqn:theta-shifts-general-k}
\begin{matrix}
    \theta(v+k\tau,\tau) = (-1)^k e^{-2\pi \im k ( v + k\tau/2)  } \theta(v,\tau)
    ,\\
    \theta_1(v+k\tau,\tau) =  e^{-2\pi \im k ( v + k\tau/2)} \theta(v,\tau)  
    ,\\
\theta_2(v+k\tau,\tau) = (-1)^k  e^{-2\pi \im k ( v + k\tau/2)} \theta_2(v,\tau)  
,\\
\theta_3(v+k\tau,\tau) =   e^{-2\pi \im k ( v + k\tau/2)} \theta_3(v,\tau)  
.
\end{matrix}
\end{equation}
The method to prove Equation \ref{eqn:theta-shifts-general-k} for e.g.\ $\theta$ for $k\geq0$ using Equation \ref{eqn:theta-shifts}
is 
\begin{align*}
    \theta(v+k\tau,\tau) 
    &=  (-1)^k q^{-k/2} e^{-2\pi \im  \sum_{r=0}^{k-1}(v+r\tau)} \theta(v,\tau)
\\&= (-1)^k q^{-k/2} e^{-2\pi \im ( kv + k(k-1)\tau/2)} \theta(v,\tau)
\\&= (-1)^k e^{-2\pi \im ( kv + k(k-1)\tau/2) +  (-k\pi \im \tau)  } \theta(v,\tau)
\\&= (-1)^k e^{-2\pi \im ( kv + k^2\tau/2)  } \theta(v,\tau)
,
\end{align*}
and so, with $k\geq0$, to find the transformation for $-k$ we see
\begin{align*}
    \theta(v,\tau) 
    &=  (-1)^k q^{-k/2} e^{-2\pi \im  \sum_{r=1}^{k}(v-r\tau)} \theta(v-k\tau,\tau)
\\&= (-1)^k q^{-k/2} e^{-2\pi \im ( kv - k(k+1)\tau/2)} \theta(v-k\tau,\tau)
\end{align*}
and similarly for $\theta_\lambda$ when $\lambda=1,2,3$ . 
From this we see the formulas above hold when $k$ is positive or nonpositive.
    
\vspace{1em}

The following lemma follows the idea of Lemma 1.1 in \cite{KL-EGTF} :

\begin{lem}    
\label{lemma:double-period}



Let $\lambda \in \{1,2,3\}$, let $(t,\tau) \in \mathbb{C}\times\mathbf{H}$ be in the domain of $\mathscr{L}_\lambda(t,\tau)$.
Then:

\begin{itemize}
    \item[(1)] $\mathscr{L}_\lambda(t,\tau) = \mathscr{L}_\lambda(t + a,\tau)$
for any $a \in 2\mathbb{Z}$.
\item[(2)] If $p_1(V)_{S^1}=0$,
   then $\mathscr{L}_\lambda(t,\tau) = \mathscr{L}_\lambda(t + a \tau ,\tau)$
for any $a \in 2\mathbb{Z}$.
\end{itemize}

\begin{proof}

(1): 
Consider the replacement $t \to t+a$ for some $a \in 2\mathbb{Z}$.
From the expression for $\mathscr{L}_\lambda$
found in Lemma \ref{lem:lef-expression},
we see that this replacement results in a
shift in the 1st argument of each theta function in Lemma \ref{lem:lef-expression}
by an even integer,
since each of the $m_{\alpha,i}$
and $n_{\alpha,j}$'s are integers. 
Since a theta function is unchanged by a shift in its 1st argument
by an even integer (which is seen from Equation \ref{eqn:theta-shifts}),
it follows that $\mathscr{L}_\lambda$ is unchanged.

(2):
It remains to show that $\mathscr{L}_\lambda$ is unchanged under $t \to t+ a\tau$ for $a$ an even integer,
provided $p_1(V)_{S^1}=0$.
Under the replacement $t \to t+a\tau$, 
from Equation \ref{eqn:theta-shifts-general-k},
we see that in the expression for $\mathscr{L}_\lambda$ in Lemma \ref{lem:lef-expression},
the factors in the product over $i=1,\ldots,d$
gain exponential factors that cancel in the numerator and denominator,
whereas the $j=1,\ldots,\ell$ product gains an overall exponential factor of
$$
\exp \left( -2\pi \im a \left( 
    \sum_{j=1}^{\ell} m_{\alpha,j} b_{\alpha,j} 
    + (t+\frac{\tau}{2}) \sum_{j=1}^{\ell} m_{\alpha,j}^2
\right) \right) .
$$

We see that the $\mathscr{L}$ in full is unchanged if, for every connected component $M_\alpha$ of the fixed locus of $M$,
the conditions
$
\sum_{j=1}^{\ell} m_{\alpha,j} b_{\alpha,j}  = 0 
\text{ and }
\sum_{j=1}^{\ell} m_{\alpha,j}^2 = 0 
$
hold. But since each $m_{\alpha,j}^2$ is nonnegative, the second condition alone is sufficient:
$
\sum_{j=1}^{\ell} m_{\alpha,j}^2 = 0 .
$

We now show that the condition $p_1(V)_{S^1}=0$ implies the condition $\sum_{j=1}^{\ell} m_{\alpha,j}^2 = 0$. Indeed,
for each $M_\alpha$ we have
\begin{equation*}
\frac{i_\alpha^* p_1(V)_{S^1}}{(2\pi \im)^2}
 = \sum_{j=1}^{\ell} \left(\left(b_{\alpha,j}+m_{\alpha,j} \frac{w}{2\pi \im} \right)^2 \right)
= \frac{i_\alpha^*p_1(V)}{(2\pi \im)^2} 
+ 2\frac{w}{2\pi\im} \sum_{j=1}^{\ell} m_{\alpha,j} b_{\alpha,j} 
+ \frac{w^2}{(2\pi \im)^2} \sum_{j=1}^{\ell} m_{\alpha,j}^2 ,
\end{equation*}
where $w \in H^2_{S^1}$ is the universal element.
So if $p_1(V)_{S^1}=0$
then $0 = \frac{i_\alpha^* p_1(V)}{(2\pi \im)^2} + 2\frac{w}{2\pi\im} \sum_{j=1}^{\ell} m_{\alpha,j} b_{\alpha,j} 
+ \frac{w^2}{(2\pi \im)^2} \sum_{j=1}^{\ell} m_{\alpha,j}^2$.
Next, since $M_\alpha$
has the trivial $S^1$-action, 
consequently $H^*_G(M_\alpha) \cong H^*(M_\alpha) \otimes H^*(BG) $
where $G=S^1$,
and so the different powers of $w$ must be linearly independent over $H^*(M_\alpha)$,
from which
it follows that $
\sum_{j=1}^{\ell} m_{\alpha,j} b_{\alpha,j} = 0 \text{ and }\sum_{j=1}^{\ell} m_{\alpha,j} ^2 = 0
$, as was to be shown.
\end{proof}

\end{lem}

\SubSection{Modular properties of $\mathscr{L}_1,\mathscr{L}_2,\mathscr{L}_3$}

Throughout, we will reference an action of $\mathrm{SL}_2(\bb{Z})$ on $\mathbb{C} \times \mathbf{H}$
given by
$$
\left.
\begin{bmatrix}
    a&b\\c&d
\end{bmatrix}
 \cdot (v,\tau)
 = \left( \frac{v}{c\tau+d} , \frac{a\tau+b}{c\tau+d} \right)
 \right._.
$$
We will refer to this as the ``modular'' action on $\mathbb{C} \times \mathbf{H}$ (even if it does not descend to an action by the ``true'' modular group $\mathrm{PSL}_2(\bb{Z})$).
We will often reference the elements $S =  \left[\begin{smallmatrix}
0 & -1 \\
1 &  0 \\
\end{smallmatrix} \right] , T=\left[\begin{smallmatrix}
1 & 1 \\
0 &  1 \\
\end{smallmatrix} \right]$ which form a standard generating set for $\mathrm{SL}_2(\bb{Z})$. 

The induced actions of $S$ and $T$ on the Jacobi theta functions are described in the following table (see e.g.\ Section V.8 of \cite{Chand}):
\begin{equation}
\label{theta-mod-transf}
\begin{matrix}
    \theta(\frac{v}{\tau},-\frac{1}{\tau}) = \frac{1}{\im}\sqrt{\frac{\tau}{\im}} e^{\im \pi v^2 / \tau} \theta(v,\tau)  
    ,& \theta(v,\tau+1)=e^{\im\pi/4}\theta(v,\tau) 
    ,\\
    \theta_1(\frac{v}{\tau},-\frac{1}{\tau}) =  \sqrt{\frac{\tau}{\im}} e^{\im \pi v^2 / \tau} \theta_2(v,\tau)  
    ,& \theta_1(v,\tau+1)=e^{\im\pi/4}\theta_1(v,\tau) 
    ,\\
    \theta_2(\frac{v}{\tau},-\frac{1}{\tau}) =  \sqrt{\frac{\tau}{\im}} e^{\im \pi v^2 / \tau} \theta_1(v,\tau)  
    ,& \theta_2(v,\tau+1)= \theta_3(v,\tau) 
    ,\\
    \theta_3(\frac{v}{\tau},-\frac{1}{\tau}) =  \sqrt{\frac{\tau}{\im}} e^{\im \pi v^2 / \tau} \theta_3(v,\tau)  
    ,& \theta_3(v,\tau+1)= \theta_2(v,\tau) 
    ,\\
    \theta'(0,-\frac{1}{\tau}) =  \left(\frac{\tau}{\im}\right)^{3/2} \theta'(0,\tau)  
    ,& \theta'(0,\tau+1)= e^{\im\pi/4} \theta'(0,\tau) 
    .
\end{matrix}
\end{equation}

Let $\mathcal{F}$ denote the set of all holomorphic complex-valued functions on $\bb{C}\times\mathbf{H}$
that vanish nowhere.

\begin{lem} \label{lemma:modular-transf-of-L}

Suppose the equivariant first Pontryagin class 
    $p_1(V)_{S^1}$
of $V$ is 0.

Then the action of $S$ 
\begin{itemize}
    \item takes $\mathscr{L}_1$ 
        to the product of $\mathscr{L}_2$ with
            an element of $\mathcal{F}$,
    \item takes $\mathscr{L}_2$
    to the product of $\mathscr{L}_1$ with
            an element of $\mathcal{F}$,
    \item and takes $\mathscr{L}_3$ to 
    to the product of $\mathscr{L}_3$ itself with
            an element of $\mathcal{F}$.
\end{itemize}

The action of $T$ 
\begin{itemize}
    \item takes $\mathscr{L}_1$ to a constant multiple of itself,
    \item takes $\mathscr{L}_2$ to a constant multiple of $\mathscr{L}_3$,
    \item and takes $\mathscr{L}_3$ to a constant multiple of $\mathscr{L}_2$.
\end{itemize}

\begin{proof}

This is seen directly from the expressions for the $\mathscr{L}_\lambda$'s
in Lemma \ref{lem:lef-expression},
and the modular transformation properties of the theta functions
described in Equation
\ref{theta-mod-transf}.

We note for example that the expression for the result of the $S$-action looks like
$$
\mathscr{L}_\lambda(S(t,\tau))
= \left(\text{some element of } \mathcal{F} \right) 
    \sum_\alpha \int_{M_\alpha} e^{ 
    (i_\alpha^*(p_1(V)_{S^1}))|_{w=t}  / (4\pi \im)} (\cdots)
$$
where $(\cdots)$ represents the expression inside the $\int_{M_\alpha}$ integral in the formula we found for $\mathscr{L}_{\lambda'}$ in Lemma \ref{lem:lef-expression},
where $\lambda' = 2,1,3$ respectively if $\lambda = 1,2,3$.
Thus if $p_1(V)_{S^1}=0$
then $\mathscr{L}_\lambda(S(t,\tau)) = \left(\text{some element of } \mathcal{F} \right) \mathscr{L}_{\lambda'}(t,\tau)$.

We note that the $S$ and $T$ actions interchange the terms in the $\sum_{\mu=1}^{3}$ sum in Lemma \ref{lem:lef-expression}, but the overall sum is preserved, possibly up to a factor of an element of $\mathcal{F}$.
Also, the extra factor arising in the $T$ action on $\theta_1$ cancels in the numerator and denominator in the $\mu=1$ term of the $\sum_{\mu=1}^{3}$ summation.
\end{proof}

\end{lem}

\SubSection{Holomorphicity of $\mathscr{L}(t,\tau)$ at real values of $t$}

For simplicity, we assume the fixed locus of the $S^1$ action on $M$ has only finitely many connected components.

The following lemma is based on the ideas of Lemma 1.3 of \cite{KL-EGTF} and the Lemma on p.\ 180 of \cite{Hirz}:

\begin{lem}
\label{lemma:no-poles-at-real-t}

For $\lambda=1,2,3$, the function $\mathscr{L}_{\lambda}(t,\tau)$
is holomorphic in $t$ at least at real values of $t$ (for any $\tau$).

\begin{proof}

${}$ 

\vspace{0.5em}

    \noindent
Step 1:
    First let us expand $\Phi$
as a power series $\sum_{n=0}^{\infty} \Phi_n q^{n/2}$,
so that each $\Phi_n$
is itself a (virtual) finite-rank bundle, equivariant for the $S^1$-action
and carrying a twisted Dirac operator $\mathscr{D}_n = \mathscr{D} \otimes \Phi_n$ on the spinor bundle of $\Phi_n$, where
$\mathscr{D}$ is the Dirac operator on $M$.

We then have 
$ \mathscr{L}(t,\tau) = \sum_{n=0}^{\infty} b_n(t)q^{n/2}
$
where $b_n(t)$ is the Lefschetz number for $\mathscr{D}_n$ on $\Phi_n$.

Next, since each $\Phi_n$ is just a finite-rank bundle with elliptic operator
$\mathscr{D}_n$,
we see the Lefschetz number $b_n(t)$
is just a linear combination of characters of finite-dimensional irreducible $S^1$ representations.
Since all such characters of $S^1$ are just monomials in $z$ where $z \in S^1$ and $S^1$ is viewed as the complex unit circle,
this tells us that $b_n(t)$
is just a finite Laurent polynomial in $z,z^{-1}$ where $z = e^{2\pi \im t}$.
In particular, the only possible singularities of $b_n$ 
occur at $z=0$.

\vspace{0.5em}

    \noindent
Step 2:
Let us view $b_n$ as functions of $z= e^{2\pi \im t}$ instead of $t$.
Now, to show the holomorphicity of $\mathscr{L}(z,q) = \sum_{n=0}^\infty b_n(z) q^{n/2}$
in $t$ at real $t$,
it suffices to show the holomorphicity in $z$
for $z \in S^1 \subset \mathbb{C}$.

By the lemma on p.\ 180 (in Appendix III) of \cite{Hirz},
since (from Step 1)
each $b_n(z)$ already has no singularities on the unit circle,
it now suffices to show for each fixed $\tau_0 \in \mathbf{H}$ that
$\mathscr{L} = \sum_{n=0}^{\infty} b_n(z) q_0^{n/2}$
(where $q_0^{1/2} = e^{\pi\im \tau_0}$)
has at most isolated singularities for $z$ on the unit circle $S^1$.

Indeed (with the simplifying assumption that there are only finitely many components $M_\alpha$ of the fixed locus $M^{S^1}$),
the function $\mathscr{L}$
is a finite sum over $\alpha \in \mathcal{I}$ of 
an integral over $M_\alpha$ of a ratio of (finite) polynomials 
in theta functions with arguments of the form $(c+nt,\tau)$
for some integer $n$ (and cohomology class $c$).
Therefore any possible singularities of $\mathscr{L}$
must be at the zeros of some polynomial in the 
 holomorphic functions $\theta_\nu(c+ n t,\tau)$, but
necessarily the set of such zeros is discrete.
\end{proof}
    
\end{lem}


\SubSection{Proof of rigidity of $\mathscr{D}$ on $\Phi_\lambda$}

The following lemma is analogous to the discussion at the end of Section 1 of \cite{KL-EGTF} :

\begin{lem}
\label{lemma:poles-move-to-real-line}

Assume the $S^1$-action on $M$ has only \emph{isolated} fixed points.
    Let $\lambda \in \{1,2,3\}$.
If $\mathscr{L}_\lambda(t,\tau)$ has a singularity at some $(t,\tau) \in \mathbb{C} \times \mathbf{H}$,
then there is some $\lambda' \in \{1,2,3\}$
such that $\mathscr{L}_{\lambda'}(t,\tau)$ has a singularity at a \emph{real} value of $t$ (for some $\tau \in \mathbf{H}$).

\begin{proof}

First, since we are assuming the action has isolated fixed points,
it follows that in the expression for $\mathscr{L}$ in Lemma 5.6,
all the $c_{\alpha,i}$ and $b_{\alpha,j}$ are just 0,
and we can drop the integral over $M_\alpha$ (since now each $M_\alpha$ is just a single point).

Suppose for some $\tau_0 \in \mathbf{H}$ that $\mathscr{L}(t,\tau_0)$ 
has a singularity in $t$ at some location $t_0 \in \mathbb{C}$.
From the expressions for $\mathscr{L}$ in the isolated-fixed-point case,
we see this can only occur at some zero of $\theta_\nu(mt,\tau)$ for some integer $m$ and some $\nu\in\{0,1,2,3\}$.

According to the discussion on p.\ 66 of \cite{Chand},
every zero of $\theta_\nu(v,\tau)$ occurs at a location where $v \in \bb{Z} + \tau \bb{Z}$.

Thus we see a pole of $\mathscr{L}(t,\tau)$ can only occur at a location where $t \in \bb{Q} + \tau \bb{Q} $.

Now suppose $\mathscr{L}$ has a singularity at some location $(t=t_0,\tau=\tau_0)$
where $t_0 = \frac{j}{k}(c\tau_0+d)$, where $j,k,c,d$ are integers with $k\neq0$ and $c,d$ are coprime (which can be achieved by ``pulling out'' all their common integer factors and ``absorbing'' these into $j$).
Since $c,d$ are coprime we can find integers $a,d$ with $ad-bc=1$.
Let $g := \left[ \begin{smallmatrix}
    a & b \\ c & d
\end{smallmatrix} \right] \in \mathrm{SL}_2(\bb{Z})$.
Then
$g\cdot (t_0,\tau_0) = (\frac{t_0}{c\tau_0+d},\frac{a\tau_0+b}{c\tau_0+d})
= (\frac{j}{k},\frac{a\tau_0+b}{c\tau_0+d})
$;
this implies $\mathscr{L}(g\cdot (t,\tau))$ has a singularity at $t = \frac{j}{k} \in \mathbb{R}$ (for some $\tau=c\tau_0+d$).

But by Lemma \ref{lemma:modular-transf-of-L} we know $\mathscr{L}_{\lambda}(g \cdot (t,\tau))$ is
just $f(t,\tau)\mathscr{L}_{\lambda'}(t,\tau)$
for some nowhere-vanishing function $f$, and some $\lambda'\in\{0,1,2,3\}$.
\end{proof}

\end{lem}

We are finally in a position to prove Theorem \ref{main-rigidity-thm}
which we state in a slightly more precise form: 


\begin{thm}


    Suppose $p_1(V)_{S^1}=0$. 
    Suppose also the $S^1$-action on $M$ has only isolated fixed points. 
    Then, for each $\lambda=1,2,3$, 
    the function $\mathscr{L}_\lambda(t,\tau)$
    is constant in $t$, for each fixed $\tau \in \mathbf{H}$.

\begin{proof}

We already know $\mathscr{L} = \mathscr{L}_\lambda(t,\tau)$
is meromorphic,
and we showed in Lemma \ref{lemma:double-period}
that $\mathscr{L}$ is doubly-periodic in $t$.
Since a standard theorem from complex analysis says
a holomorphic doubly-periodic function on the complex plane is constant,
it only remains to show that $\mathscr{L}$ has no poles.

Indeed, if $\mathscr{L}$
has a pole, then by Lemma \ref{lemma:poles-move-to-real-line},
$\mathscr{L}$ has a pole at some real value of $t$ (for some $\tau$);
but this contradicts Lemma \ref{lemma:no-poles-at-real-t} which states that $\mathscr{L}$ has no poles at real values of $t$ (for any $\tau$).
Thus the theorem is proved.
\end{proof}

\end{thm}

\Section{Further comments}

A potential further direction to pursue
would be to try to find a generalization
of Theorem 1 of \cite{KL-Mod}
to provide another proof of Theorem \ref{main-rigidity-thm} of this paper.
Using the notation of \cite{KL-Mod}, an appropriate generalization of that paper's Theorem 1 would need to allow the consideration of $\psi(E,V)$'s
with multiple different $V$'s, in this case both $V$ and $TM$.
Towards this end, one possible methodology 
could be to investigate 
further the behavior
of the specific $L\mathrm{Spin}$
characters involved in the analysis in that paper.

\section*{Acknowledgments}

 I would like to thank my advisors, Professors Ko Honda and Kefeng Liu,
 for their support and guidance throughout the writing of this paper.
 In addition, I thank Professors Sucharit Sarkar and Peter Petersen for their help as members of my doctoral committee.

\end{document}